\documentclass[10pt]{amsart}

\usepackage{enumerate}
\usepackage{graphics}
\usepackage{graphicx}
\usepackage{amsmath}
\usepackage{amssymb}
\usepackage{amstext}
\usepackage{times}
\usepackage[latin1]{inputenc}
\usepackage{latexsym}
\usepackage{amsthm}
\usepackage{amscd}

\def\R{\mathbb{R}}
\def\L{\mathbb{L}}
\def\C{\mathbb{C}}
\def\N{\mathbb{N}}

\def\H{\mathbb{H}}

\newcommand{\Real}{\mbox{\rm Re}}
\newcommand{\Ima}{\mbox{\rm Im}}

\newcommand{\dist}{\operatorname{dist}}
\newcommand{\length}{\operatorname{length}}

\newcommand{\intc}{\operatorname{Int}}
\newcommand{\escpro}[1]{\langle {#1} \rangle}
\newcommand{\norma}[1]{\| {#1} \|}
\newcommand{\escproR}[1]{\langle {#1} \rangle_{0}}
\newcommand{\normaR}[1]{\| {#1} \|_{0}}

\def\a{{\alpha}}
\def\t{{\theta}}
\def\g{{\gamma}}
\def\l{{\lambda}}
\def\de{{\delta}}

\def\Om{{\Omega}}
\def\s{{\sigma}}
\def\ep{{\epsilon}}
\def\ri{{\rm i}}

\newtheorem{lemma}{Lemma}
\newtheorem{remark}{Remark}
\newtheorem{theorem}{Theorem}

\newtheorem{definition}{Definition}
\newtheorem{claim}{Claim}[section]

\numberwithin{equation}{section}

\begin{document}

\title{On the existence of a proper conformal maximal disk in $\L^3$}
\author{Antonio Alarcón}

\thanks{This research is partially supported by MEC-FEDER Grant no. MTM2004 - 00160.}
\date{\today}

\address{Departamento de Geometría y Topología\hfill\break\indent  Universidad de Granada, \hfill\break\indent18071, Granada \hfill\break\indent Spain}

\email{alarcon@ugr.es}

\begin{abstract}
In this paper we construct an example of a properly immersed maximal surface in the Lorentz-Minkowski space $\L^3$ with the conformal type of a disk.
\\


\end{abstract}

\maketitle


\section{Introduction}

The conformal type problem has strongly influenced the modern theory of surfaces (see for instance \cite{Perez}, \cite{Weitsman}, \cite{cast}, \cite{kus}), in particular, the theory of maximal surfaces in the Lorentz-Minkowski space $\L^3$ (see \cite{Isa-PL2}, \cite{Calabi}, \cite{Um-Ya}, \cite{Isa-PL1}, among others). This paper is closely related to an intrinsic question associated with the underlying complex structure, the type problem for a maximal surface, i.e., to determine whether its conformal structure is parabolic or hyperbolic. The family of all open Riemann surfaces can be divided into three mutually exclusive classes: elliptic (i.e. compact), parabolic and hyperbolic. A Riemann surface without boundary is called hyperbolic if it carries a non-constant positive superharmonic function and parabolic if it is neither compact nor hyperbolic (see \cite{a-s} or \cite{Farkas} for details).

A maximal hypersurface in a Lorentzian manifold is a spacelike hypersurface with zero mean curvature. Besides of their mathematical interest these hypersurfaces and more generally those having constant mean curvature have a significant importance in physics (cf. \cite{Kiehn}, \cite{Kiehn2}, \cite{MT}). When the ambient space is the Minkowski space $\L^n,$
one of the most important results is the proof of a Bernstein-type theorem for maximal hypersurfaces
in $\L^n$. Calabi proved in \cite{Calabi} that the only complete hypersurfaces with zero mean curvature in $\L^3$ (i.e. maximal surfaces) and $\L^4$ are spacelike hyperplanes, solving the so called Bernstein-type problem in dimensions 3 and 4. Cheng and Yau in \cite{Ch-Yau} extended this result to $\L^n,$ $n \geq 5.$
It is therefore meaningless to consider global problems on maximal and everywhere regular hypersurfaces in $\L^n.$ In contrast, there exists a lot of results about existence of non-flat parabolic maximal surfaces with singularities (see for example \cite{Isa-PL1}, \cite{Isa-PL2}, \cite{Isa-PL-RS}).
In this paper, we construct the first example of a proper maximal surface in $\L^3$ with singularities and with hyperbolic type. We would like to point out that our example does not have branch points, all the singularities are of lightlike type (see definition \ref{def:light} in page \pageref{pag:def}).

More precisely, we prove the following existence theorem.

\begin{theorem}\label{teorema}
There exists a conformal proper maximal immersion of the disk (with lightlike singularities).
\end{theorem}

For several reasons, lightlike singularities of maximal surfaces in $\L^3$ are specially interesting. This kind of singularities are more interesting than branch points, in the sense that they have a physical interpretation (see \cite{Kiehn}, \cite{Kiehn2}). At these points, the limit tangent plane is lightlike, the curvature blows up and the Gauss map has no well defined limit. However, if we allow branch points, then proving the analogous result of Theorem \ref{teorema} has less technical difficulties.
\\

The fundamental tools used in the proof of this result (Runge's theorem and the López-Ros transformation) were previously utilized by Morales in \cite{Morales} to construct the first example of a proper minimal surface in $\R^3$ with the conformal type of a disk. This technique was improved by Martín and Morales \cite{MM-Convex}, \cite{MM-Convex2} and later by Alarcón, Ferrer and Martín \cite{AFM} in order to construct hyperbolic minimal surfaces in $\R^3$ with more complicated topology.
\\

{\em Acknowledgments.} We would like to thank F. J. López for helpful criticisms of the paper.


\section{Background and notation}

\subsection{The Lorentz-Minkowski space}

We denote by $\L^3$ the three dimensional Lorentz-Minkowski space $(\R^3,\escpro{\cdot,\cdot}),$ where $\escpro{\cdot,\cdot}=dx_1^2+dx_2^2-dx_3^2.$ The Lorentzian {\em norm} is given by $\|(x_1,x_2,x_3)\|^2=x_1^2+x_2^2-x_3^2.$ We say that a vector $v\in\R^3\setminus \{(0,0,0)\}$ is spacelike, timelike or lightlike if $\|v\|^2$ is positive, negative or zero, respectively. The vector $(0,0,0)$ is spacelike by definition. A plane in $\L^3$ is spacelike, timelike or lightlike if the induced metric is Riemannian, non degenerate and indefinite or degenerate, respectively.

In order to differentiate between $\L^3$ and $\R^3$, we denote $\R^3=(\R^3,\escproR{\cdot,\cdot}),$ where $\escproR{\cdot,\cdot}$ is the usual metric of $\R^3$, i.e., $\escproR{\cdot,\cdot}=dx_1^2+dx_2^2+dx_3^2.$ We also denote the Euclidean norm by $\normaR{\cdot}.$
\\

By an (ordered) $\L^3$-orthonormal basis we mean a basis of $\R^3,$ $\{u,v,w\},$ satisfying
\begin{itemize}
\item $\escpro{u,v}=\escpro{u,w}=\escpro{v,w}=0$;
\item $\norma{u}=\norma{v}=-\norma{w}=1$.
\end{itemize}
Notice that $u$ and $v$ are spacelike vectors whereas $w$ is timelike. In addition, we say that an $\L^3$-orthonormal basis is peculiar if $\escproR{u,v}=\escproR{v,w}=0.$ In particular, $\{u,v,w\}$ is a peculiar $\L^3$-orthonormal basis if and only if $\{u,v,w\}$ is an $\L^3$-orthonormal basis and the third coordinate of $v$ is zero. In that case, we also have $\normaR{v}=1.$
\\

We call $\H^2:=\{(x_1,x_2,x_3)\in\R^3 \;|\; x_1^2+x_2^2-x_3^2=-1\}$ the hyperbolic sphere in $\L^3$ of constant intrinsic curvature $-1.$ Notice that $\H^2$ has two connected components $\H^2_+:=\H^2\cap\{x_3\geq 1\}$ and $\H^2_-:=\H^2\cap \{x_3\leq -1\}.$ The stereographic projection $\s$ for $\H^2$ is the map $\s:\H^2\to \C\cup\{\infty\}\setminus \{|z|=1\}$ given by
\[
\s (x_1,x_2,x_3)=\frac{x_1+\ri x_2}{1-x_3} \;,\quad \s(0,0,1)=\infty\;.
\]

\subsection{Translating spheres}

Given a real number $r,$ we define
\[
b(r):= (0,0,r)+\H^2_-=\{(x_1,x_2,x_3)\in\R^3 \;|\; \norma{(x_1,x_2,x_3-r)}^2=-1,\; x_3\leq r-1\}\;.
\]
We also define
\begin{multline*}
B(r):=(0,0,r)+\{(x_1,x_2,x_3)\in\R^3 \;|\; \norma{(x_1,x_2,x_3)}^2<-1,\; x_3\leq -1\}=
\\
\{(x_1,x_2,x_3)\in\R^3 \;|\; \norma{(x_1,x_2,x_3-r)}^2<-1,\; x_3\leq r-1\}\;.
\end{multline*}
Notice that $b(0)=\H^2_-$ and $b(r)=\partial B(r).$ Moreover, if $r_1<r_2,$ then $\overline{B(r_1)}\subset B(r_2)$ and $b(r_1)\cap b(r_2)=\emptyset.$ Furthermore, $\R^3=\cup_{r\in\R}B(r).$
\\

Now, for $r\in\R,$ consider the set
\[
E(r):= \{(x_1,x_2,x_3)\in\R^3 \;|\; x_1^2+x_2^2>0,\; x_3 \leq r-1\}
\]
Notice that there is a bijection $ [0,+\infty[\times ]0,+\infty[\times [-\pi,\pi[\rightarrow E(r)$ given by
\[
(s,t,\t)\longmapsto (t \cos \t, t\sin\t, r-\sqrt{s^2+1})\;.
\]
Observe that $\overline{B(r)}\setminus \{x_1^2+x_2^2=0\}$ is included in $E(r).$ At this point, we define the horizontal projection to $b(r)$ as the map $\mathcal{P}_H^r:E(r)\to b(r)$ given by
\[
\mathcal{P}_H^r (t \cos \t, t\sin\t, r-\sqrt{s^2+1})= (s \cos \t, s\sin\t, r-\sqrt{s^2+1})\;.
\]
Observe that this map does not depend on $t.$ Using this projection, we can define another two maps which aim to out $B(r).$ First, we define $\mathcal{N}_H^r:E(r)\to \{(x_1,x_2,0)\in\R^3 \;|\; x_1^2+x_2^2=1\}$ as
\[
\mathcal{N}_H^r (t \cos \t, t\sin\t, r-\sqrt{s^2+1}) = (\cos \t, \sin \t, 0)\;.
\]
Notice that $\mathcal{N}_H^r$ neither depends on $t$ nor $s$ and
\[
\mathcal{N}_H^r (p) = \frac{\mathcal{P}_H^r(p)-p}{\norma{\mathcal{P}_H^r(p)-p}} = \frac{\mathcal{P}_H^r(p)-p}{\normaR{\mathcal{P}_H^r(p)-p}}\;, \quad \forall p\in B(r)\cap \{x_1^2+x_2^2>0\}\;.
\]
Note that, for any $p\in E(r),$ one has that $\mathcal{N}_H^r(p)$ is a spacelike vector with $\norma{\mathcal{N}_H^r (p)}=\normaR{\mathcal{N}_H^r (p)}=1.$

Consider now $\mathcal{N}^r:b(r)\to \H^2_+$ the exterior $\L^3$-normal Gauss map of $b(r).$ Then, we define the map $\mathcal{N}_N^r:E(r)\to \H^2_+$ as
\[
\mathcal{N}_N^r(p) = \mathcal{N}^r(\mathcal{P}_H^r(p))\;.
\]
Notice that $\mathcal{N}_N^r(p)$ is a timelike vector for all $p\in E(r)$ and
\[
\mathcal{N}_N^r(t \cos \t, t \sin \t, r-\sqrt{s^2+1}) = (-s \cos \t, -s \sin \t, \sqrt{s^2+1})\;,
\]
therefore, $\mathcal{N}_N^r$ does not depend on $t.$
\begin{figure}[htbp]
    \begin{center}
        \includegraphics[width=0.80\textwidth]{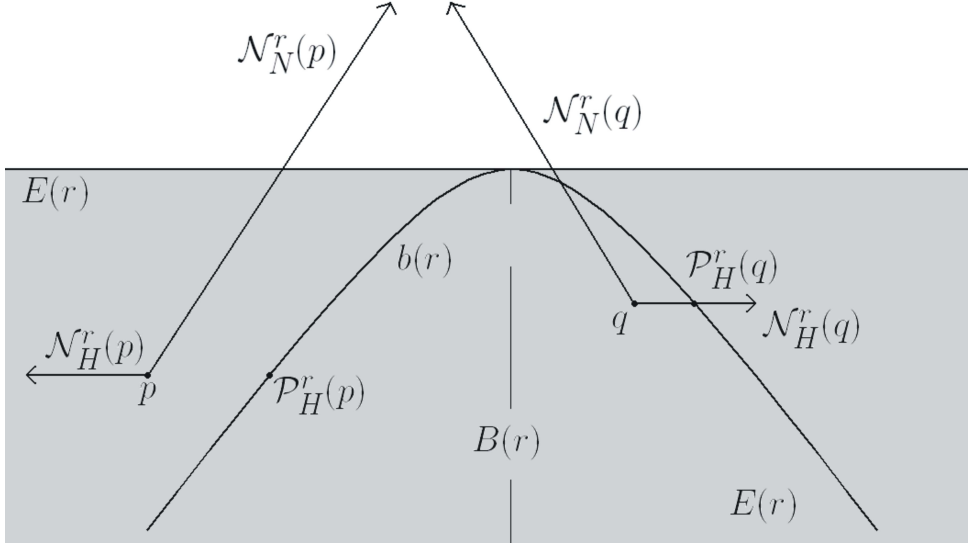}
    \end{center}
    \caption{The set $E(r)$ and the two associated maps.} \label{fig:normales}
\end{figure}

\subsection{Maximal surfaces}

Any conformal maximal immersion $X:M\to\L^3$ is given by a triple $\Phi=(\Phi_1,\Phi_2,\Phi_3)$ of holomorphic 1-forms defined on the Riemann surface $M,$ having no common zeros and satisfying
\begin{equation}\label{2}
|\Phi_1|^2+|\Phi_2|^2-|\Phi_3|^2\neq 0\;;
\end{equation}
\begin{equation}\label{conforme}
\Phi_1^2+\Phi_2^2-\Phi_3^2=0\;;
\end{equation}
and all periods of the $\Phi_j$ are purely imaginary. Here we consider $\Phi_i$ to be a holomorphic function times $dz$ in a local parameter $z.$ Then, the maximal immersion $X:M\to\L^3$ can be parameterized by $z\mapsto \Real \int^z \Phi.$ The above triple is called the Weierstrass representation of the maximal immersion $X.$ Usually, the second requirement \eqref{conforme} is guaranteed by the introduction of the formulas
\[
\Phi_1=\frac{\ri}{2}(1-g^2)\eta\;,\quad \Phi_2=-\frac12(1+g^2)\eta\;,\quad \Phi_3=g\eta
\]
for a meromorphic function $g$ with $|g(p)|\neq 1,$ $\forall p\in M,$ (the stereographically projected Gauss map) and a holomorphic 1-form $\eta.$ We also call $(g,\eta)$ or $(g,\Phi_3)$ the Weierstrass representation of $X.$

In this paper, we deal with maximal immersions with lightlike singularities, according with the following definition.
\begin{definition}\label{def:light} \label{pag:def}
A point $p\in M$ is a lightlike singularity of the immersion $X$ if $|g(p)|=1.$
\end{definition}

In this article, all the maximal immersions are defined on simply connected domains of $\C,$ thus the Weierstrass 1-forms have no periods and so the only requirements are \eqref{2} at the points that are not singularities, and \eqref{conforme}. In this case, the differential $\eta$ can be written as $\eta =f(z)dz.$ The metric of $X$ can be expressed as
\begin{equation}\label{metrica}
ds^2=\frac12 (|\Phi_1|^2+|\Phi_2|^2-|\Phi_3|^2)=\big( \frac12 (1-|g|^2)|f||dz| \big)^2\;.
\end{equation}
The Euclidean metric on $\C$ is denoted as $\escpro{,}=|dz|^2.$ Note that $ds^2=\l_X^2\,|dz|^2,$ where the conformal coefficient $\l_X$ is given by \eqref{metrica}.

Along this paper, we use some $\L^3$-orthonormal bases. Given $X:\Omega\to\L^3$ a maximal immersion and $S$ an $\L^3$-orthonormal basis, we write the Weierstrass data of $X$ in the basis $S$ as
\[
\Phi_{(X,S)}=(\Phi_{(1,S)},\Phi_{(2,S)},\Phi_{(3,S)})\;, \quad f_{(X,S)}\;, \quad g_{(X,S)}\;, \quad \eta_{(X,S)}\;.
\]

In the same way, given $v\in\R^3,$ we denote by $v_{(k,S)}$ the $k$th coordinate of $v$ in $S.$ We also represent by $v_{(*,S)}=(v_{(1,S)},v_{(2,S)})$ the first two coordinates of $v$ in the basis $S.$

Given a curve $\a$ in $\Om,$ by $\length (\a,ds)$ we mean the length of $\a$ with respect to the metric $ds.$ Given a subset $W\subset \Om,$ we define
\begin{itemize}
\item $\dist_{(W,ds)}(p,q)=\inf \{\length (\a,ds) \;|\; \a:[0,1]\to W,$ $\a(0)=p,$ $\a(1)=q\},$ for any $p,q\in W.$

\item $\dist_{(W,ds)}(U,V)= \inf \{ \dist_{(W,ds)}(p,q) \;|\; p\in U,$ $q\in V\},$ for any $U,V\subset W.$
\end{itemize}

\subsection{The López-Ros transformation}

The proof of Lemma \ref{lema} exploits what has come to be call the López-Ros transformation. If $(g,f)$ are the Weierstrass data of a maximal immersion $X:\Omega\to\L^3$ (being $\Omega$ simply connected), we define on $\Omega$ the data
\[
\widetilde{g}=\frac{g}{h}\;,\quad \widetilde{f}=f\,h\;,
\]
where $h:\Omega\to\C$ is a holomorphic function without zeros. Observe that the new meromorphic data satisfy \eqref{2} at the regular points, and \eqref{conforme}, so the new data define a maximal immersion (possibly with different lightlike singularities) $\widetilde{X}:\Omega\to\L^3.$ This method provides us with a powerful and natural tool for deforming maximal surfaces. One of the most interesting properties of the resulting surface is that the third coordinate function is preserved.


\section{Proof of Theorem \ref{teorema}}\label{se:teo}

In order to prove Theorem \ref{teorema} we will apply the following technical Lemma. It will be proved later in Section \ref{se:lema}.

\begin{lemma}\label{lema}
Let $r_1$ and $r_2$ be two real numbers with $r_1<r_2.$ Let $\mathcal{O}\subset\C$ be a simply connected domain such that $0\in\mathcal{O},$ and consider $X:\mathcal{O}\to\L^3$ a non-flat conformal maximal immersion (possibly with lightlike singularities) with $X(0)=0.$ Suppose that there exists a polygon $P\subset\mathcal{O}$ satisfying $0\in\intc P$ and
\begin{equation}\label{6}
X(\mathcal{O}\setminus \intc P)\subset B(r_2)\setminus \overline{B(r_1)}\;.
\end{equation}
Consider $b_1>0.$ Then, for any $b_2>0$ such that $r_2-b_2>r_1,$ there exist a polygon $Q$ and a non-flat conformal maximal immersion (possibly with lightlike singularities) $Y:\overline{\intc Q}\to\L^3$ satisfying:
\begin{enumerate}[\rm (I)]
\item $P\subset \intc Q\subset \overline{\intc Q}\subset \mathcal{O}$.
\item $Y(0)=0$.
\item $\|Y(z)-X(z)\|_{0}<b_1,$ $\forall z\in \overline{\intc P}$.
\item $Y(Q)\subset B(r_2)\setminus \overline{B(r_2-b_2)}$.
\item $Y(\intc Q\setminus \intc P)\subset \L^3\setminus B(r_1-1-b_2)$.
\end{enumerate}
\end{lemma}

Using this Lemma, we will construct a sequence of immersions $\{\psi_n\}_{n\in\N}$ that converges to an immersion $\psi$ which proves Theorem \ref{teorema}, up to a reparametrization of its domain. Previously, consider $\{s_n\}$ a sequence of real numbers given by
\[
s_1>1\;,\quad s_k=s_{k-1}+\frac{2}{k},\; k>1\;.
\]
Notice that this sequence diverges. We also consider another sequence of reals $\{\a_n\}$ satisfying
\[
\prod_{k=1}^\infty \a_k =\frac12\;,\quad 0<\a_k<1\;,\quad \forall k\in\N \;.
\]

Now, we are going to construct a sequence $\{\Upsilon_n\}_{n\in\N},$ where the element $\Upsilon_n=\{U_n,\psi_n,P_n\}$ consists of an open domain $U_n,$ a non-flat conformal maximal immersion $\psi_n:U_n\to\L^3$ and a polygon $P_n$ on $\C.$

We construct the sequence in order to satisfy the following properties:
\begin{enumerate}[\rm (A{$_n$})]
\item $0\in \intc P_n \subset \overline{\intc P_n}\subset U_n.$

\item $\psi_n(0)=(0,0,0).$

\item $\overline{\intc P_{n-1}}\subset \intc P_n \subset \overline{\intc P_n} \subset \mathcal{D},$ where $\mathcal{D}$ is a bounded simply connected domain of $\C$ which does not depend on $n.$

\item $\psi_n(P_n)\subset B(s_n)\setminus \overline{B(s_n-\frac1{(n+1)^2})}.$

\item $\psi_n(\intc P_n\setminus \intc P_{n-1}) \subset \L^3 \setminus B(s_{n-1}-\frac1{n^2}-1-\frac1{(n+1)^2}).$

\item $\normaR{\psi_n(z)-\psi_{n-1}(z)}<\frac1{n^2},$ $\forall z\in \overline{\intc P_{n-1}}.$

\item $\l_{\psi_n}(z)\geq \a_n\l_{\psi_{n-1}}(z),$ $\forall z\in \overline{\intc P_{n-1}}.$
\end{enumerate}

The sequence $\{\Upsilon_n\}$ is constructed in a recursive way. The existence of a non-flat conformal maximal immersion $\psi_1:U_1\to \L^3$ and a polygon $P_1$ satisfying (A$_1$), (B$_1$) and (D$_1$) is straightforward. The rest of the properties have no sense for $n=1.$

Assume we have got $\Upsilon_1,\ldots,\Upsilon_{n-1}.$ We are going to construct $\Upsilon_n.$ We choose a decreasing sequence of positive reals $\{\ep_m\}_{m\in\N}\searrow 0$ with $\ep_m<1/n^2$ for all $m\in\N.$ For each $m,$ we consider the immersion $Y_m$ and the polygon $Q_m$ given by Lemma \ref{lema} for the following data:
\[
X=\psi_{n-1}\,,\quad P=P_{n-1}\;,\quad r_1=s_{n-1}-\frac1{n^2}\;,\quad r_2=s_n\;,\quad b_1=\ep_m\;,\quad b_2=\frac1{(n+1)^2}\;,
\]
and $\mathcal{O}$ a simply connected domain with $\overline{\intc P_{n-1}}\subset \mathcal{O} \subset U_{n-1} \subset\mathcal{D}$ and satisfying \eqref{6}. The existence of this domain is a consequence of (D$_{n-1}$). From (III) in Lemma \ref{lema}, we deduce that the sequence $\{Y_m\}$ uniformly converges to $\psi_{n-1}$ on $\overline{\intc P_{n-1}}.$ Then, taking into account that $Y_m$ is a harmonic map and that its metric is given by its derivatives, we conclude that the sequence $\{\l_{Y_m}\}$ uniformly converges to $\l_{\psi_{n-1}}$ on $\overline{\intc P_{n-1}}.$ Hence, there exists $m_0\in \N$ satisfying
\begin{equation}\label{lan}
\l_{Y_{m_0}}(z)\geq \a_n \l_{\psi_{n-1}}(z)\;,\quad \forall z\in \overline{\intc P_{n-1}}\;.
\end{equation}
Then, we define $\psi_n:=Y_{m_0}$ and $P_n:=Q_{m_0}.$ Properties (A$_n$) and (C$_n$) are consequence of (I) whereas (B$_n$), (D$_n$), (E$_n$) and (F$_n$) are obtained from (II), (IV), (V) and (III), respectively. Finally, \eqref{lan} implies (G$_n$). This concludes the construction of the sequence $\{\Upsilon_n\}.$
\\

Now, define $\Delta:=\cup_{n\in\N} \intc P_n.$ Since (C$_n$), $\Delta$ is a bounded simply connected domain of $\C,$ i.e., $\Delta$ is biholomorphic to a disk. Moreover, from (F$_n$) we obtain that $\{\psi_n\}$ is a Cauchy sequence, uniformly on compact sets of $\Delta.$ Then, Harnack's Theorem guarantees the existence of a harmonic map $\psi:\Delta\to\L^3$ such that $\{\psi_n\}\to\psi,$ uniformly on compact sets of $\Delta.$ Moreover, $\psi$ has the following properties:
\\

{\em $\bullet$ $\psi$ is maximal and conformal.}
\\

{\em $\bullet$ $\psi$ is an immersion:} Indeed, for any $z\in\Delta$ there exists $n_0\in \N$ so that $z\in \intc P_{n_0}.$ Given $k>n_0$ and using (G$_j$), $j=n_0+1,\ldots,k,$ one has
\[
\l_{\psi_k}(z)\geq \a_k\cdots \a_{n_0+1}\l_{\psi_{n_0}}(z)\geq \a_k\cdots \a_1\l_{\psi_{n_0}}(z)\;.
\]
Taking the limit as $k\to\infty,$ we infer that
\[
\l_\psi (z)\geq \frac12 \l_{\psi_{n_0}}(z)>0\;,
\]
and so, $\psi$ is an immersion.
\\

{\em $\bullet$ $\psi$ is proper in $\L^3$:} Consider $K\subset\L^3$ a compact set. For each $n\in\N,$ define
\[
t_n:= s_{n-1}-\frac1{n^2}-1-\frac1{(n+1)^2}\;.
\]
Notice that $t_n>s_{n-1}-3,$ and so $\{t_n\}$ diverges. Then, for any positive constant $\xi,$ there exists $n_0\in\N$ satisfying
\[
K\subset B(t_n-\xi)\;,\quad \forall n\geq n_0\;.
\]
From properties (E$_n$), we have
\begin{equation}\label{reina}
\psi_n(z)\in \L^3\setminus B(t_n)\;,\quad \forall z\in \intc P_n\setminus \intc P_{n-1}\;.
\end{equation}
If we fix a large enough $\xi>0,$ and taking (F$_k$), $k\geq n,$ into account, we obtain from \eqref{reina} that
\[
\psi(z)\in \L^3 \setminus B(t_n-\xi) \;,\quad \forall z\in \intc P_n\setminus \intc P_{n-1}\;.
\]
Then we have $\psi^{-1}(K)\cap (\intc P_n\setminus \intc P_{n-1})=\emptyset,$ for $n\geq n_0.$ Therefore, $\psi^{-1}(K)\subset \intc P_{n_0-1},$ and so it is compact in $\Delta.$
\\

This completes the proof of Theorem \ref{teorema}.

\begin{remark}
$\psi$ is proper in $\L^3$ and it has the conformal type of a disk. Therefore, $\psi$ is a non-flat immersion.
\end{remark}


\section{Proof of Lemma \ref{lema}}\label{se:lema}

Throughout the proof, we will use the following two constants:
\begin{itemize}
\item $\mu=\sup\{\dist_{\R^3}(p,\mathcal{P}_H^{r_2}(p)) \;|\; p\in b(r_1)\cap E(r_2)\}=\sqrt{(r_2-r_1)^2+2(r_2-r_1)}.$ Notice that since $r_1<r_2$ we have that $b(r_1)\cap E(r_2) = b(r_1)\setminus \{(0,0,r_1-1)\}.$

\item $\ep_0>0$ is taken small enough to satisfy all the inequalities appearing in this section. This choice depends only on the data of the lemma.
\end{itemize}

\subsection{Preparing the first inductive process}

\begin{claim}\label{Af-1}
There exists a simply connected domain $W,$ with $\overline{\intc P}\subset W\subset \overline{W}\subset \mathcal{O}\subset \C,$ and there exists a set of points $\{p_1,\ldots,p_n\}\subset (W\setminus \overline{\intc P})\cap E(r_2)$ (for some $n\in\N$) satisfying the following list of properties:
\begin{enumerate}[\rm i)]
\item Labeling $p_{n+1}=p_1,$ the segments $\overline{p_1p_2},\ldots,\overline{p_{n-1}p_n},\overline{p_np_{n+1}}$ form a polygon $\widehat{P}\subset W\setminus \overline{\intc P}$.

\item For each $i\in\{1,\ldots,n\},$ there exists an open disk $B^i\subset W\setminus \overline{\intc P}$ such that $\{p_i,p_{i+1}\}\subset B^i$ and
\begin{equation}\label{7}
\normaR{X(z)-X(w)}<\ep_0\;,\quad \forall z,w\in B^i\;.
\end{equation}

\item For each $i\in\{1,\ldots,n\},$ there exists a peculiar $\L^3$-orthonormal basis $S_i=\{e^i_1,e^i_2,e^i_3\}$ with
\[
e^i_2=\mathcal{N}_H^{r_2}(X(p_i))\;, \quad e^i_3=(0,0,1)\;,
\]
and satisfying
\begin{equation}\label{8}
\normaR{e^i_j-e^{i+1}_j}<\frac{\ep_0}{3\mu}\;,\quad \forall j\in\{1,2,3\}\;,
\end{equation}
and
\begin{equation}\label{9}
f_{(X,S_i)}(p_i)\neq 0\;.
\end{equation}

\item For each $i\in\{1,\ldots,n\},$ there exists a complex number $\t_i$ such that $|\t_i|=1,$ $\Ima (\t_i)\neq 0$ and
\begin{equation}\label{10'}
\left| \overline{\t_i\, \frac{f_{(X,S_i)}(p_i)}{|f_{(X,S_i)}(p_i)|}}+1\right|< \frac{\ep_0}{3\mu}\;.
\end{equation}
\end{enumerate}
\end{claim}
\begin{proof}
If the points $p_i,$ $i=1,\ldots,n$ are taken close enough and the natural number $n$ is sufficiently large, then the existence of the simply connected domain $W$ and properties i), ii) and \eqref{8} are a direct consequence of the uniform continuity of $X$ and $\mathcal{N}_H^{r_2}$ in $E(r_2)\setminus L,$ where $L$ is a small open neighborhood of $\{x_1^2+x_2^2=0\}.$

Now, observe that if $f_{(X,S_i)}(p_i)=0,$ then $g_{(X,S_i)}(p_i)=\infty.$ Since $X$ is non-flat, this fact only occurs in a finite set of points. Therefore, we can choose the points satisfying \eqref{9}.

Finally, the choice of the complex numbers $\t_i$ satisfying iv) is straightforward.
\end{proof}

\begin{remark}\label{R3-ort}
Observe that for any $i\in\{1,\ldots,n\}$ the peculiar $\L^3$-orthonormal basis $S_i$ is also an orthonormal basis of $\R^3,$ i.e.,
\begin{equation}\label{base-ort}
\escproR{e^i_j,e^i_k}=0\;,\forall j\neq k\quad \text{and}\quad \normaR{e^i_j}=1\;,\forall j=1,2,3\;.
\end{equation}
\end{remark}

The proof of the following claim is straightforward.

\begin{claim}\label{Af-a}
There exists $\de\in ]0,1[$ small enough satisfying the following properties:
\begin{enumerate}[\rm ({a}1)]
\item $\overline{\intc \widehat{P}\setminus \cup_{k=1}^n D(p_k,\de)}$ is a simply connected domain, where we denote by $D(p_k,\de)$ the disk centered at $p_k$ with radius $\de.$

\item $\overline{D(p_i,\de)\cup D(p_{i+1},\de)}\subset B^i,$ $\forall i=1,\ldots,n.$

\item $\overline{D(p_i,\de)}\cap \overline{D(p_k,\de)}=\emptyset,$ $\forall \{i,k\}\subset \{1,\ldots,n\},$ $i\neq k.$

\item $\de\cdot \max_{\overline{D(p_i,\de)}} \{|f_{(X,S_i)}|\}<2\ep_0,$ $\forall i=1,\ldots,n.$

\item $\de\cdot \max_{\overline{D(p_i,\de)}} \{|f_{(X,S_i)} g^2_{(X,S_i)}|\}< 2\ep_0 |\Ima\, \t_i|,$ $\forall i=1,\ldots,n.$

\item $\de\cdot \max_{\overline{D(p_i,\de)}} \{\normaR{\phi}\}<\ep_0,$ $\forall i=1,\ldots,n,$ where $\Phi=\phi\, dz$ is the Weierstrass representation of the maximal immersion $X$.

\item $3\mu\, \max_{w\in \overline{D(p_i,\de)}} \{|f_{(X,S_i)}(w)-f_{(X,S_i)}(p_i)|\}<\ep_0\, |f_{(X,S_i)}(p_i)|,$ $\forall i=1,\ldots,n.$
\end{enumerate}
\end{claim}

Now, define
\begin{equation}\label{12}
\ell := \sup\{ \dist_{(\mathcal{E},\escpro{,})}(0,z) \;|\; z\in\mathcal{E}\}+2\pi\de+\de+1\;,
\end{equation}
where
\[
\mathcal{E}:=\overline{\intc \widehat{P}\setminus\cup_{k=1}^n D(p_k,\de)}\;.
\]

\subsection{The first inductive process}

The first inductive process consists of the construction of a sequence $\Psi_1,\ldots,\Psi_n,$ where the element $\Psi_i=\{k_i,a_i,C_i,G_i,\Phi^i,D_i\}$ is composed of the following ingredients:
\begin{itemize}
\item $k_i$ is a suitable positive real constant.

\item $a_i$ is a point lying on the segment $\overline{p_i q_i},$ where $q_i=p_i+\de.$ Notice that $\overline{p_i q_i}\setminus \{q_i\}\subset D(p_i,\de)$ and $q_i\in \partial D(p_i,\de).$

\item $C_i$ is an arc of the circumference centered at $p_i$ that contains the point $a_i.$

\item $G_i$ is a closed annular sector bounded by $C_i,$ a piece of $\partial D(p_i,\de)$ and two radii of this circumference.

\item $\Phi^i$ is a Weierstrass representation on $\overline{W}.$ We also write $\Phi^i=\phi^i\, dz,$ where $\phi^i:\overline{W}\to\C^3$ is a meromorphic map. The points $p_1,\ldots,p_i$ will be poles of $\Phi^i.$

\item $D_i$ is a simply connected domain of $\C$ satisfying $\overline{D_i}\cap G_i=\emptyset$ and
\[
\{p_i,w_i:=p_i-k_i\t_i\}\subset D_i\subset \overline{D_i}\subset D(p_i,\de)\;.
\]
\end{itemize}
\begin{figure}[htbp]
    \begin{center}
        \includegraphics[width=0.5\textwidth]{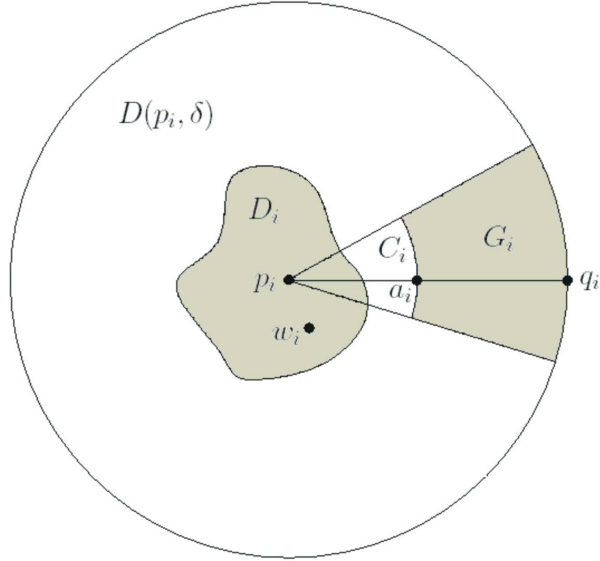}
    \end{center}
    \caption{The disk $D(p_i,\de).$} \label{fig:disco}
\end{figure}

\begin{remark}
From now on, we will use the convention that $\Psi_{n+1}=\Psi_1.$
\end{remark}

\begin{claim}\label{Af-b}
We can construct the sequence satisfying the following properties:
\begin{enumerate}[\rm ({b}1.{$i$})]
\item $\de\cdot \max_{\overline{D(p_k,\de)}} \{|f_{(\Phi^i,S_k)}|\}<2\ep_0,$ $\forall k=i+1,\ldots,n.$

\item $\de\cdot \max_{\overline{D(p_k,\de)}} \{ |f_{(\Phi^i,S_k)}\, g^2_{(\Phi^i,S_k)}| \}<2\ep_0\cdot |\Ima\,\t_k|,$ $\forall k=i+1,\ldots,n.$

\item $3\mu\cdot \max_{w\in\overline{D(p_k,\de)}} \{ |f_{(\Phi^i,S_k)}(w)-f_{(X,S_k)}(p_k)| \}<\ep_0\cdot |f_{(X,S_k)}(p_k)|,$ $\forall k=i+1,\ldots,n.$

\item $\normaR{\Real \int_{\a_z}\Phi^i}<\ep_0,$ $\forall z\in C_i,$ where $\a_z$ is a piece of $C_i$ joining $a_i$ and $z.$

\item $\Phi^i_{(3,S_i)}=\Phi^{i-1}_{(3,S_i)}.$

\item $\normaR{\phi^i(z)-\phi^{i-1}(z)}<\frac{\ep_0}{n\,\ell},$ $\forall z\in \overline{W}\setminus \big( D(p_i,\de)\cup (\cup_{k=1}^{i-1}D_k) \big).$

\item $|\Real\int_{\overline{q_i z}}\Phi^i_{(1,S_i)}|<4\ep_0,$ $\forall z\in G_i.$

\item $\big| \Real\int_{\overline{q_i z}}\Phi^i_{(2,S_i)}-\frac12 \big( \int_{\overline{q_i z}} \frac{k_i\, dw}{w-p_i}\big) |f_{(X,S_i)}(p_i)| \big|<4\ep_0,$ $\forall z\in G_i.$

\item $\normaR{\Real\int_{\overline{q_i a_i}} \Phi^i -\Real\int_{\overline{q_{i-1} a_{i-1}}}\Phi^{i-1}}<21\ep_0.$
\end{enumerate}
\end{claim}

All the above properties have meaning for $i=1,\ldots,n$ except (b1.$i$), (b2.$i$) and (b3.$i$), which hold only for $i=1,\ldots,n-1.$ Similarly, (b9.$i$) only occurs for $i=2,\ldots,n+1.$ Notice that properties (b5.$i$), (b7.$i$) and (b8.$i$) tell us that the deformation of our surface around the points $p_i$ follows the direction of $e^i_2=\mathcal{N}_H^{r_2}(X(p_i)).$

As we have announced, we construct the family $\Psi_1,\ldots, \Psi_n$ in a recursive way. Let $\Phi^0=\phi\, dz$ be the Weierstrass representation of the immersion $X.$ We denote $\Psi_0=\{\Phi^0\}.$ Suppose we have constructed $\Psi_1,\ldots,\Psi_{i-1}.$ We are going to construct $\Psi_i.$

The Weierstrass data $\Phi^i$, in the basis $S_i$, are determined by the López-Ros transformation
\[
f_{(\Phi^i,S_i)}= f_{(\Phi^{i-1},S_i)}\cdot h_i\;,\quad g_{(\Phi^i,S_i)}=\frac{g_{(\Phi^{i-1},S_i)}}{h_i}\;,
\]
where $h_i:\overline{W}\to \overline{\C}$ is given by
\[
h_i(z)=\frac{k_i\,\t_i}{z-p_i}+1\;.
\]

We choose the constant $k_i>0$ small enough to satisfy properties (b1.$i$), (b2.$i$), (b3.$i$) and (b6.$i$). Notice that this choice is possible since $\Phi^i$ converges uniformly to $\Phi^{i-1}$ on $\overline{W}\setminus (D(p_i,\de)\cup (\cup_{k=1}^{i-1} D_k))$ if $k_i\to 0,$ and since we can use (b1.$i-1$), (b2.$i-1$) and (b3.$i-1$). In the case $i=1,$ these properties are consequence of (a4), (a5) and (a7).

\begin{remark}
The meromorphic function $h_i$ is close to 1 outside a neighborhood of $p_i.$ The constant $\t_i$ has the effect of a rotation near to $p_i.$ This rotation let us to choose the direction of deformation of the surface. Outside a neighborhood of $p_i$ this effect almost disappears.
\end{remark}

Furthermore, property (b5.$i$) trivially follows from the definition of $\Phi^i.$

We choose $a_i$ as the first point in the (oriented) segment $\overline{q_i p_i}$ that satisfy
\begin{equation}\label{13}
\frac12 |f_{(X,S_i)}(p_i)| \int_{\overline{q_i a_i}} \frac{k_i\, dw}{w-p_i}=3\mu\;.
\end{equation}

Let $D_i$ be a simply connected domain containing the pole, $p_i,$ and the zero, $w_i=p_i-k_i\t_i,$ of $h_i$ and satisfying $\overline{D_i}\subset D(p_i,\de)$ and $\overline{D_i}\cap \overline{q_i a_i}=\emptyset.$ We can take it because $w_i\notin \overline{p_i q_i}$ (recall that $\Ima\, \t_i \neq 0$).

Before proving (b7.$i$) and (b8.$i$) we are going to check the following inequality:
\begin{equation}\label{b78}
\left\| \left( \Real\int_{\overline{q_i z}} \Phi^i_1 \right) e^i_1 + \left( \Real\int_{\overline{q_i z}} \Phi^i_2\right) e^i_2\right.
\left. -\frac12 \left( \int_{\overline{q_i z}} \frac{k_i\, dw}{w-p_i}\right) |f_{(X,S_i)}(p_i)|\,e^i_2\right\|_{0}<4\ep_0\;,\quad \forall z\in \overline{q_i a_i}\;.
\end{equation}
Consider $z\in \overline{q_i a_i}.$ Taking \eqref{13} and \eqref{10'} into account, we obtain
\begin{equation}\label{jessica}
\left| \frac12 \left( \int_{\overline{q_i z}} \frac{k_i\, dw}{w-p_i}\right) |f_{(X,S_i)}(p_i)|+ \frac12 \left( \int_{\overline{q_i z}} \frac{k_i\, dw}{w-p_i}\right) \overline{\t_i\, f_{(X,S_i)}(p_i)}\right|<
3\mu \left| 1+ \frac{\overline{\t_i\, f_{(X,S_i)}(p_i)}}{|f_{(X,S_i)}(p_i)|}\right|<\ep_0\;.
\end{equation}
For convenience we use complex notation and we write $a+\ri b$ instead of $a e^i_1+b e^i_2$ (recall that $S_i$ is an $\R^3$-orthonormal basis). Then, taking into account \eqref{jessica} and the fact that $\Real\, \Phi_1+\ri \Real\,\Phi_2=-\frac{\ri}2(f\,g^2\,dw+\overline{f}\,d\overline{w}),$ we obtain that
\[
\left| \left( \Real\int_{\overline{q_i z}} \Phi^i_{(1,S_i)} \right) + \ri \left( \Real\int_{\overline{q_i z}} \Phi^i_{(2,S_i)} \right)- \frac{\ri}2 \left( \int_{\overline{q_i z}} \frac{k_i\, dw}{w-p_i}\right) |f_{(X,S_i)}(p_i)| \right|<
\]
\[
\frac12 \left| \int_{\overline{q_i z}} \overline{f_{(\Phi^i,S_i)}(w)\, dw} + \int_{\overline{q_i z}} f_{(\Phi^i,S_i)}(w)\, g^2_{(\Phi^i,S_i)}(w)\, dw \right.
\left.-\left( \int_{\overline{q_i z}} \frac{k_i\, dw}{w-p_i}\right) \overline{\t_i\, f_{(X,S_i)}(p_i)} \right|+\ep_0\;.
\]
Using the definition of $\Phi^i$ and $h_i,$ the last expression is less than
\begin{multline*}
\frac12 \left| \int_{\overline{q_i z}} \overline{\big( f_{(\Phi^{i-1},S_i)}(w)-f_{(X,S_i)}(p_i) \big)\frac{k_i\,\t_i}{w-p_i}\,dw}\right| + \frac12 \left|\int_{\overline{q_i z}} \overline{f_{(\Phi^{i-1},S_i)}(w)\, dw}\right|
\\
+\frac12 \left| \int_{\overline{q_i z}} f_{(\Phi^{i-1},S_i)}(w)\, g^2_{(\Phi^{i-1},S_i)}(w)\,\frac{dw}{h_i(w)}\right|+\ep_0<4\ep_0\;,
\end{multline*}
where we have used \eqref{13}, (b1.$i-1$), (b2.$i-1$), (b3.$i-1$) and the fact that $|h_i(w)|>|\Ima\, \t_i|,$ $\forall w\in \overline{q_i a_i}.$ Thus, we have proved that \eqref{b78} holds. Therefore, if $C_i$ and $G_i$ are chosen sufficiently close to $a_i$ and $\overline{q_i a_i}$, respectively, we obtain \eqref{b78} for all $z\in G_i$ and (b4.$i$). Now, Properties (b7.$i$) and (b8.$i$) follow straightforwardly.

Finally, in order to prove (b9.$i$) we write
\[
\left\| \Real \int_{\overline{q_i a_i}} \Phi^i - \Real \int_{\overline{q_{i-1} a_{i-1}}} \Phi^{i-1} \right\|_{0}\leq
\sum_{j=1}^3 \left\| \left( \Real \int_{\overline{q_i a_i}} \Phi^i_{(j,S_i)}\right) e^i_j - \left( \Real \int_{\overline{q_{i-1} a_{i-1}}} \Phi^{i-1}_{(j,S_{i-1})} \right) e^{i-1}_j \right\|_{0}\;,
\]
and we are going to bound each addend separately. Using that $S_i$ is an orthonormal basis of $\R^3,$ (b7.$i$) and (b7.$i-1$), we have
\begin{multline*}
\left\| \left( \Real \int_{\overline{q_i a_i}} \Phi^i_{(1,S_i)}\right) e^i_1 - \left( \Real \int_{\overline{q_{i-1} a_{i-1}}} \Phi^{i-1}_{(1,S_{i-1})} \right) e^{i-1}_1 \right\|_{0} <
\\
 \left| \Real \int_{\overline{q_i a_i}} \Phi^i_{(1,S_i)} \right| + \left| \Real \int_{\overline{q_{i-1} a_{i-1}}} \Phi^{i-1}_{(1,S_{i-1})} \right|<4\ep_0+4\ep_0=8\ep_0\;.
\end{multline*}
For $j=2,$ we use (b8.$i$), (b8.$i-1$), \eqref{13} and \eqref{8} to obtain
\[
\left\| \left( \Real \int_{\overline{q_i a_i}} \Phi^i_{(2,S_i)}\right) e^i_2 - \left( \Real \int_{\overline{q_{i-1} a_{i-1}}} \Phi^{i-1}_{(2,S_{i-1})} \right) e^{i-1}_2 \right\|_{0} <
\]
\[
\left\| \frac12 \left( |f_{(X,S_i)}(p_i)|\int_{\overline{q_i a_i}} \frac{k_i\, dw}{w-p_i} \right) e^i_2 - \frac12 \left( |f_{(X,S_{i-1})}(p_{i-1})|\int_{\overline{q_{i-1} a_{i-1}}} \frac{k_{i-1}\, dw}{w-p_{i-1}} \right) e^{i-1}_2 \right\|_{0} +8\ep_0 =
\]
\[
\left\| 3\mu e^i_2 - 3\mu e^{i-1}_2 \right\|_{0}+8\ep_0 < \ep_0+8\ep_0 =9\ep_0\;.
\]
For the last addend, we use (b5.$i$), (b5.$i-1$) and the fact that $\normaR{e^i_3}=\normaR{e^{i-1}_3}=1$ to obtain
\[
\left\| \left( \Real \int_{\overline{q_i a_i}} \Phi^i_{(3,S_i)}\right) e^i_3 - \left( \Real \int_{\overline{q_{i-1} a_{i-1}}} \Phi^{i-1}_{(3,S_{i-1})} \right) e^{i-1}_3 \right\|_{0} \leq
\]
\[
\left| \Real \int_{\overline{q_i a_i}} \Phi^i_{(3,S_i)} \right| + \left| \Real \int_{\overline{q_{i-1} a_{i-1}}} \Phi^{i-1}_{(3,S_{i-1})} \right|<
\]
\[
\de\cdot \big( \max_{\overline{D(p_i,\de)}} \big\{\normaR{\phi^{i-1}} \big\} + \max_{\overline{D(p_{i-1},\de)}} \big\{ \normaR{\phi^{i-2}} \big\} \big) <\de \big( \frac{\ep_0}{\ell} + \frac{\ep_0}{\de} + \frac{\ep_0}{\ell} + \frac{\ep_0}{\de} \big) < 4\ep_0\;,
\]
where we have used (b6.$k$) for $k=1,\ldots,i-1$ and (a6).

\subsection{Preparing the second inductive process}

Note that the Weierstrass representations $\Phi^i$ have simples poles and zeros in $W.$ Our next step consists of describing a simply connected domain $\Om$ in $W$ where the above Weierstrass representations define maximal immersions with lightlike singularities.

For each $i\in\{1,\ldots,n\},$ consider $D^i$ an open disk centered at $p_i$ and so that $D_i\subset \overline{D(p_i,\de)}\subset D^i$ and $D^1,\ldots,D^n$ are pairwise disjoint. Let $\a_i$ be a  simple curve in $D^i\setminus \overline{D(p_i,\de)}$ connecting $q_i$ with $\partial D^i\cap \intc \widehat{P}.$ Finally, take $N_i$ a small open neighborhood of $\a_i\cup\overline{q_i a_i}$ included in $\overline{G_i \cup (D^i\setminus D(p_i,\de))}.$

\begin{claim}\label{Af-C}
If $D^i,$ $\a_i$ and $N_i$ are suitably chosen, then the domain $\Om$ given by
\[
\Om:=\big( \intc \widehat{P}\setminus \cup_{k=1}^n D^k \big) \cup \big( \cup_{k=1}^n N_k \big)\;,
\]
satisfies the following properties:
\begin{enumerate}[\rm ({c}1)]
\item $\overline{\Om}$ is a simply connected domain.

\item $\overline{q_i a_i}\subset \overline{\Om}$ and $\overline{\intc P}\subset \Om.$

\item $\overline{\Om}$ does not contain any pole $p_i$ and any zero $w_i$ of the function $h_i,$ for all $i=1,\ldots,n.$

\item $\sup_{z\in\overline{\Om}}\{\dist_{(\Om,\escpro{,})}(0,z)\}<\ell,$ where $\ell$ has been defined on \eqref{12}.

\item $\overline{\Om}\cap \overline{D(p_i,\de)}\subset G_i.$
\end{enumerate}
\end{claim}

Taking (c1) and (c3) into account, we can define $n$ maximal immersions with lightlike singularities $X_1,\ldots,X_n,$
where $X_i:\Om'\to\L^3$ is given by
\[
X_i(z)=\Real \int_0^z \phi^i(w)\, dw\;,
\]
where $\Om'$ is a suitable open neighborhood of $\overline{\Om}$ satisfying (c4).

\begin{claim}\label{Af-D}
For $i=1,\ldots,n,$ we have
\begin{enumerate}[\rm ({d}1.${i}$)]
\item $\normaR{X_i(z)-X_{i-1}(z)}<\ep_0/n,$ $\forall z\in \Om'\setminus D(p_i,\de).$

\item $(X_i)_{(3,S_i)}=(X_{i-1})_{(3,S_i)}.$

\item $\normaR{X_n(a_i)-X_n(a_{i+1})}<26\ep_0.$

\item $\normaR{X_n(a_i)-(X(p_i)+3\mu\, \mathcal{N}_H^{r_2}(X(p_i)))}<14\ep_0.$

\item $X_n(a_i)\in\L^3\setminus B(r_1+2(r_2-r_1)).$
\end{enumerate}
\end{claim}
\begin{proof}
In order to get (d1.$i$) we use (b6.$i$) and (c4) as follows:
\[
\normaR{X_i(z)-X_{i-1}(z)}=\left\| \Real \int_0^z(\phi^i-\phi^{i-1})\,dz \right\|_{0} \leq
\left| \int_0^z \left\| \phi^i-\phi^{i-1} \right\|_{0} \,dz \right|\leq \frac{\ep_0}{n\ell} \left|\int_0^z dz\right|\leq \frac{\ep_0}{n}\;.
\]

(d2.$i$) is a direct consequence of (b5.$i$). In order to check (d3.$i$), we apply (d1.$k$), $k=1,\ldots,n,$ \eqref{7} and (b9.$i+1$) to obtain
\begin{multline*}
\normaR{X_n(a_i)-X_n(a_{i+1})}\leq \normaR{X_n(a_i)-X_i(a_i)} + \normaR{X_n(a_{i+1})-X_{i+1}(a_{i+1})}
\\
 + \normaR{X_{i+1}(q_{i+1})-X(q_{i+1})}+\normaR{X(q_{i+1})-X(q_i)}+\normaR{X(q_i)-X_i(q_i)}
\\
+\normaR{(X_i(a_i)-X_i(q_i))-(X_{i+1}(a_{i+1})-X_{i+1}(q_{i+1}))}<
\end{multline*}
\[
4\ep_0 + \normaR{X(q_{i+1})-X(q_i)}+ \left\| \Real\int_{\overline{q_{i+1} a_{i+1}}} \Phi^{i+1}-\Real\int_{\overline{q_i a_i}}\Phi^i \right\|_{0}< 4\ep_0+\ep_0+21\ep_0=26\ep_0\;.
\]

Now, we are going to prove (d4.$i$). Using (d1.$k$), $k=i+1,\ldots,n,$ one gets
\[
\normaR{X_n(a_i)-(X(p_i)+3\mu\, \mathcal{N}_H^{r_2}(X(p_i)))}\leq
\]
\[
\normaR{X_n(a_i)-X_i(a_i)} + \normaR{X_i(a_i)-(X_i(q_i)+3\mu\,\mathcal{N}_H^{r_2}(X(p_i)))} + \normaR{X_i(q_i)-X(p_i)}<
\]
\begin{multline*}
\ep_0+\normaR{(X_i(a_i)-X_i(q_i))_{(*,S_i)}-3\mu\, e^i_2} + |(X_i(a_i)-X_i(q_i))_{(3,S_i)}|
\\
+\normaR{X_i(q_i)-X(q_i)}+\normaR{X(q_i)-X(p_i)}<\ep_0+8\ep_0+2\ep_0+2\ep_0+\ep_0=14\ep_0\;,
\end{multline*}
where we have used (b7.$i$), (b8.$i$), (b5.$i$) and \eqref{7}.

Finally, (d5.$i$) is a consequence of (d4.$i$) and the fact that $X(p_i)+3\mu\, \mathcal{N}_H^{r_2}(X(p_i))$ belongs to $\L^3\setminus B(r_1+3(r_2-r_1)).$
\end{proof}

In the second inductive process, we employ new basis. For each $i=1,\ldots,n,$ we take $T_i=\{w^i_1,w^i_2,w^i_3\}$ a peculiar $\L^3$-orthonormal basis so that
\begin{equation}\label{15}
w^i_3=\mathcal{N}_N^{r_2}(X_n(a_i))\;.
\end{equation}

\begin{remark}\label{loca}
Notice that, if $\ep_0$ is small enough, then $w^i_3$ is well defined, i.e. $X_n(a_i)\in E(r_2),$ because of (d4.$i$) and \eqref{6}. Observe that $\{w^i_1,w^i_2\}$ is a basis of the tangent plane $\Pi_i$ to $b(r_2)$ at $\mathcal{P}_H^{r_2}(X_n(a_i)).$ Since $B(r_2)$ is convex, we have that $B(r_2)$ is contained in a connected component of $\L^3\setminus \Pi_i,$ even more, $(p-\mathcal{P}_H^{r_2}(X_n(a_i)))_{(3,T_i)}\leq 0,$ $\forall p\in \overline{B(r_2)}.$
\end{remark}

Given $i\in \{1,\ldots,n\},$ we define $Q_i$ as the connected component of $\overline{\partial \Om\setminus (C_i \cup C_{i+1})}$ that does not cut $C_k$ for all $k\notin\{i,i+1\}.$ Observe that $\{Q_i\;|\;i=1,\ldots,n\}$ satisfy $\overline{Q_i}\cap\overline{Q_j} =\emptyset,$ for all $i\neq j,$ and the following properties:
\begin{equation}\label{16}
Q_i\subset B^i\;\text{(recall that we defined $B^i$ in Claim \ref{Af-1})}\;,
\end{equation}
\begin{equation}\label{17}
Q_i\cap \overline{D(p_k,\de)}=\emptyset\;, \quad\forall k\neq\{i,i+1\}\;,
\end{equation}
and, up to a small perturbation of the curve $Q_i$,
\begin{equation}\label{18}
f_{(X_n,T_i)}(z)\neq 0\;, \quad \forall z\in Q_i\;.
\end{equation}
\begin{figure}[htbp]
    \begin{center}
        \includegraphics[width=1\textwidth]{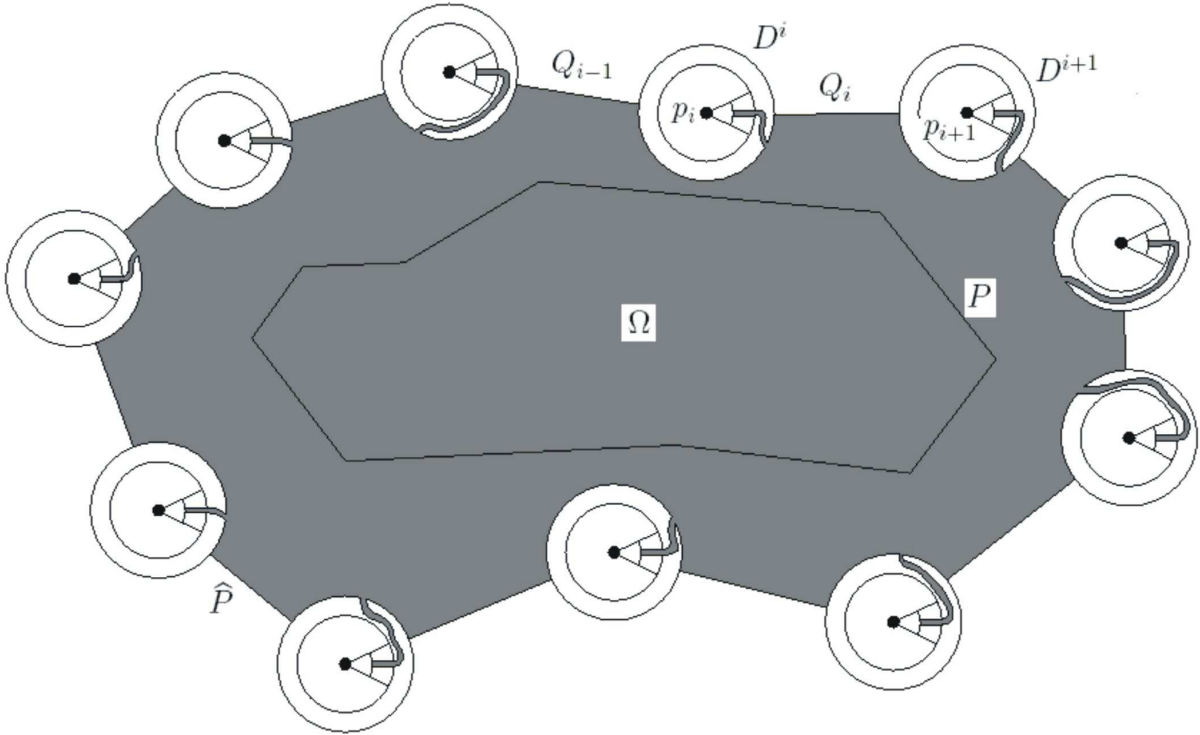}
    \end{center}
    \caption{The domain $\Om$ and the curves $Q_i.$} \label{fig:dominio}
\end{figure}

Now, for each $i=1,\ldots,n,$ let $\widehat{C}_i$ be an open set containing $C_i$ and so that
\begin{equation}\label{19}
\normaR{X_n(z)-X_n(a_i)}<3\ep_0\;,\quad \forall z\in \widehat{C_i}\cap \overline{\Om}\;.
\end{equation}
The existence of such sets is due to properties (d1.$k$), $k=i+1,\ldots,n,$ and (b4.$i$). We also define, for each $i=1,\ldots,n$ and for any $\xi>0,$ $Q_i^\xi:=\{z\in\C \;|\; \dist_{(\C,\escpro{,})}(z,Q_i)\leq \xi\}.$

\begin{claim}\label{Af-E}
There exists $\xi>0$ small enough so that:
\begin{enumerate}[\rm ({e}1)]
\item $Q_i^\xi\subset \Om'.$

\item $Q_i^\xi\cap Q_j^\xi=\emptyset,$ $\forall i\neq j.$

\item $Q_i^\xi\cap \overline{D(p_k,\de)}=\emptyset,$ $\forall k\notin \{i,i+1\}.$

\item $Q_i^\xi\subset B^i.$

\item $Q_i^{\xi/2}$ and $\overline{\Om\setminus Q_i^\xi}$ are simply connected.

\item $|f_{(X_n,T_i)}(z)-f_{(X_n,T_i)}(x)|<\ep_1,$ $\forall x\in B(z,\xi/2),$ $\forall z\in Q_i,$ where $\ep_1:=\tfrac14 \min_{Q_i}\{|f_{(X_n,T_i)}|\}.$

\item $\sup \{\dist_{(\overline{\Om\setminus Q_i^\xi},\escpro{,})} (0,z) \;|\; z\in \overline{\Om\setminus Q_i^\xi}\}<\ell.$

\item $\normaR{X_n(z)-X_n(x)}<\ep_0,$ $\forall x\in B(z,\xi/2),$ $\forall z\in Q_i.$
\end{enumerate}
\end{claim}
Observe that properties (e3), (e4) and (e7) are consequences of \eqref{17}, \eqref{16} and (c4), respectively. It is straightforward to check the other ones for a sufficiently small $\xi>0.$
\\

For each $i=1,\ldots,n,$ the plane $\Pi_i$ generated by $w^i_1$ and $w^i_2$ is spacelike. Therefore, given $z\in Q_i,$ and $v\in\Pi_i$ with $\normaR{v}=1,$ there exists $\l(v,z)\geq 0$ minimum so that
\begin{equation}\label{biblio}
X_n(z)+u+\l\cdot v\in \L^3\setminus \overline{B(r_2)}\;, \quad \forall u\in \{x_1^2+x_2^2+x_3^2\leq 1\}\;,\quad \forall \l>\l(v,z)\;.
\end{equation}
Now, we define $\Lambda:=\max \{\Lambda_1,\ldots,\Lambda_n\},$ where
\[
\Lambda_i:=\max \{\l(v,z)\;|\; z\in Q_i,\;v\in\Pi_i,\;\normaR{v}=1\}\;.
\]
Therefore, for any $u\in \L^3$ with $\normaR{u}\leq 1,$ for any $\l>\Lambda,$ and for any $i=1,\ldots,n,$ since \eqref{biblio} we obtain that
\begin{equation}\label{martilla}
X_n(z)+u+\l\cdot v\in \L^3\setminus \overline{B(r_2)}\;,\quad \forall z\in Q_i\;,\; \forall v\in \Pi_i\text{ with }\normaR{v}=1\;\;.
\end{equation}

\subsection{The second inductive process}

We are now ready to construct a sequence $\{\Xi_i \;|\; i=1,\ldots,n\},$ where the element $\Xi_i=\{Y_i,\tau_i,\nu_i\}$ is composed of:
\begin{itemize}
\item $Y_i:\Om'\to\L^3$ is a conformal maximal immersion with $Y_i(0)=0.$ We also define $Y_0:=X_n.$

\item $\{(\tau_i,\nu_i)\in \R^+\times\R^+ \;|\; i=1,\ldots,n\}.$
\end{itemize}

\begin{claim}\label{Af-F}
We can construct the sequence $\{\Xi_i \;|\; i=1,\ldots,n\}$ satisfying the following list of properties:
\begin{enumerate}[\rm ({f}1.${i}$)]
\item $(Y_i)_{(3,T_i)}=(Y_{i-1})_{(3,T_i)}.$

\item $\normaR{Y_i(z)-Y_{i-1}(z)}<\ep_0/n,$ $\forall z\in \overline{\Om\setminus Q_i^\xi}.$

\item $|f_{(Y_i,T_k)}(z)-f_{(Y_{i-1},T_k)}(z)|<\ep_1/n,$ $\forall z\in \overline{\Om\setminus Q_i^\xi},$ $\forall k=i+1,\ldots,n.$

\item $(\frac{1}{\tau_i}+\frac{\nu_i}{\tau_i(\tau_i-\nu_i)}) \max_{Q_i^\xi} \{ |f_{(Y_{i-1},T_i)}\, g^2_{(Y_{i-1},T_i)}| \} + \nu_i\, \max_{Q_i^\xi} \{ |f_{(Y_{i-1},T_i)}| \}<\frac{2}{\xi}.$

\item $\frac12(\frac{\tau_i \xi}4 \min_{Q_i} \{ |f_{(Y_0,T_i)}| \}-1)>2(\Lambda+1).$
\end{enumerate}
\end{claim}
The sequence is constructed by recursion. Consider $\Xi_0=\{Y_0\}.$ All of the properties make no sense for $i=0.$ Assume we have defined $Y_0,\ldots,Y_{i-1}.$ Then, we use Runge's theorem to get a holomorphic function without zeros, $l_i:\C\to\C,$ satisfying
\begin{itemize}
\item $|l_i(z)-\tau_i|<\nu_i,$ $\forall z\in Q_i^{\xi/2}.$

\item $|l_i(z)-1|<\nu_i,$ $\forall z \in \overline{\Om\setminus Q_i^\xi}.$
\end{itemize}
Hence, we define $Y_i(z)=\Real\int_0^z\Phi$ as the maximal immersion whose Weierstrass data in the $\L^3$-orthonormal basis $T_i$ are given by
\[
f_{(Y_i,T_i)}=f_{(Y_{i-1},T_i)}\cdot l_i\;,\quad g_{(Y_i,T_i)}=\frac{g_{(Y_{i-1},T_i)}}{l_i}\;.
\]
Note that $Y_i:\Om'\to\L^3$ is obtained from $Y_{i-1}$ by applying a López-Ros transformation. So, property (f1.$i$) trivially holds. The fact that $\phi_{(Y_i,T_k)} \stackrel{\nu_i \to 0}{\longrightarrow} \phi_{(Y_{i-1},T_k)}$ uniformly on $\overline{\Om\setminus Q_i^\xi},$ implies the rest of the properties if the constant $\nu_i$ is sufficiently small and $\tau_i$ is large enough.

\subsection{The immersion $Y$ solving Lemma \ref{lema}}

Consider the maximal immersion $Y:\Om\to\L^3$ given by $Y=Y_n.$ We are going to check that $Y$ satisfies the statements of Lemma \ref{lema}.
\\

{\bf Item (II):} It is obvious from the definition of $Y$.
\\

{\bf $Y$ is non-flat and Item (III):} Items i) and ii) in Claim \ref{Af-1} and properties (e4) and (a2) imply that
\[
\overline{\intc P}\subset \Om\setminus \left( \left( \cup_{k=1}^n D(p_k,\de) \right) \cup \left( \cup_{k=1}^n Q_k^\xi \right) \right)\;,
\]
therefore, we can successively apply (f2.$k$) and (d1.$k$), $k=1\ldots,n,$ to obtain $\forall z\in\overline{\intc P}$
\begin{equation}\label{20}
\normaR{Y(z)-X(z)}\leq \normaR{Y_n(z)-Y_0(z)}+\normaR{X_n(z)-X_0(z)}<2\ep_0<b_1\;,
\end{equation}
that proves Item (III). If $\ep_0$ is small enough, then $Y$ is non-flat because of \eqref{20} and the fact that $X$ is non-flat.
\\

{\bf Items (I) and (IV):} As a previous step we will prove the following claim:

\begin{claim}\label{Af-2}
Every connected curve $\g$ in $\Om$ connecting $P$ with $\partial \Om$ contains a point $z'\in\g$ such that $Y(z')\in \L^3\setminus\overline{B(r_2)}.$
\end{claim}
\begin{proof}
Consider $\g\subset \overline{\Om}$ a connected curve with $\g(0)\in P$ and $\g(1)=z_0\in \partial\Omega.$
\\

\noindent{\bf Case 1)} Assume $z_0\in \widehat{C}_i\cap Q_i^\xi.$ Taking Remark \ref{loca} into account, we finish proving that $(Y_n(z_0)-\mathcal{P}_H^{r_2}(X_n(a_i)))_{(3,T_i)}>0.$

Using (f2.$k$), $k\neq i,$ (f1.$i$) and \eqref{19}, we obtain
\[
|(Y_n(z_0)-X_n(a_i))_{(3,T_i)}|<4\ep_0\;.
\]
On the other hand, taking (d5.$i$) into account we know that
\[
(X_n(a_i)-\mathcal{P}_H^{r_2}(X_n(a_i)))_{(3,T_i)}>\varsigma\;,
\]
where $\varsigma$ is a positive constant depending on $r_1$ and $r_2.$ Therefore, for a small enough $\ep_0$, one has
\begin{equation}\label{21}
(Y_n(z_0)-\mathcal{P}_H^{r_2}(X_n(a_i)))_{(3,T_i)}>(X_n(a_i)-\mathcal{P}_H^{r_2}(X_n(a_i)))_{(3,T_i)}-4\ep_0>\varsigma-4\ep_0>0\;.
\end{equation}

\noindent{\bf Case 2)} Suppose $z_0\in \widehat{C}_i\cap Q_{i-1}^\xi.$ Reasoning as in the above case and using property (d3.$i-1$), one has
\begin{multline*}
|(Y_n(z_0)-X_n(a_{i-1}))_{(3,T_{i-1})}|\leq
\\
|(Y_n(z_0)-Y_0(a_i))_{(3,T_{i-1})}|+ \normaR{X_n(a_i)-X_n(a_{i-1})}< 4\ep_0+26\ep_0=30\ep_0\;.
\end{multline*}
Then, we conclude the proof in this case following the arguments of \eqref{21}.
\\

\noindent{\bf Case 3)} Assume $z_0\in \widehat{C}_i\setminus \cup_{k=1}^n Q_k.$ Since (f2.$k$), $k=1,\ldots,n$ and \eqref{19}, we obtain
\[
\normaR{Y_n(z_0)-X_n(a_i)}<4\ep_0\;,
\]
and so, if $\ep_0$ is small enough, we can finish using (d5.$i$).
\\

\noindent{\bf Case 4)} Finally, assume $z_0\in Q_i\setminus \cup_{k=1}^n C_k.$ This is the most complicated case. For the sake of simplicity, we will write $f^{i-1}$ and $g^{i-1}$ instead of $f_{(Y_{i-1},T_i)}$ and $g_{(Y_{i-1},T_i)},$ respectively. As $T_i$ is a peculiar $\L^3$-orthonormal basis, we do not lose generality using complex notation, i.e., we will write $a\eta+\ri b$ instead of $a w^i_1+b w^i_2,$ where $\eta=\normaR{w^i_1}\geq 1.$

Consider $z_1\in\g\cap \partial D(z_0,\xi/2).$ Hence, taking into account (f2.$k$), $k=i+1,\ldots,n,$ and that $\eta\geq 1,$ we have
\[
\normaR{(Y_n(z_0)-Y_n(z_1))_{(*,T_i)}}\geq \normaR{(Y_i(z_0)-Y_i(z_1))_{(*,T_i)}}-2\ep_0 =
\]
\[
\left| \left( \Real \int_{\overline{z_1 z_0}} \Phi^i_{(1,T_i)} \right) \eta +\ri \left( \Real \int_{\overline{z_1 z_0}} \Phi^i_{(2,T_i)} \right) \right| -2\ep_0 \geq
\left| \Real \int_{\overline{z_1 z_0}} \Phi^i_{(1,T_i)}  +\ri\,  \Real \int_{\overline{z_1 z_0}} \Phi^i_{(2,T_i)} \right| -2\ep_0 =
\]
using the definition of $Y_i$ and that $\Real\, \Phi_1+\ri\,\Real\Phi_2=-\frac{\ri}{2}(\overline{f}+f g^2),$ the above equation continuous
\[
\frac12 \left| \int_{\overline{z_1 z_0}} \overline{f^{i-1}l_i\, dz} + \int_{\overline{z_1 z_0}} \frac{f^{i-1}(g^{i-1})^2}{l_i}\, dz \right| -2\ep_0\geq
\]
\begin{multline*}
\frac{\tau_i}2 \left| \int_{\overline{z_1 z_0}} \overline{f^{i-1}\, dz} \right| - \frac12 \left| \frac1{\tau_i} \int_{\overline{z_1 z_0}} f^{i-1}(g^{i-1})^2\, dz \right|
\\
- \frac12 \left| \int_{\overline{z_1 z_0}} \overline{f^{i-1}(l_i-\tau_i)\, dz} \right| - \frac12 \left| \int_{\overline{z_1 z_0}} f^{i-1} (g^{i-1})^2\left( \frac1{l_i}-\frac1{\tau_i} \right)\, dz \right|-2\ep_0\geq
\end{multline*}
taking into account the definition of $l_i$ and the fact that $|\int_{\overline{z_1 z_0}}dz|=\xi/2,$
\begin{multline}\label{puti}
\frac{\tau_i}2 \left| \int_{\overline{z_1 z_0}} f^{i-1}\, dz \right| - \frac{\xi}4 \left( \frac1{\tau_i} \max_{Q_i^\xi} \{|f^{i-1}(g^{i-1})^2|\} + \nu_i\,\max_{Q_i^\xi} \{|f^{i-1}|\} \right.
\\
\left. + \frac{\nu_i}{\tau_i(\tau_i-\nu_i)}\max_{Q_i^\xi} \{|f^{i-1} (g^{i-1})^2|\}\right) -2\ep_0\geq \frac12 \left( \tau_i \left| \int_{\overline{z_1 z_0}} f^{i-1}\, dz \right| -1\right)-2\ep_0\;,
\end{multline}
where we have used (f4.$i$) in the last inequality. On the other hand, taking into account (e6), (f3.$k$), $k=1,\ldots,i-1,$ and the definition of $\ep_1,$ we can deduce
\begin{multline*}
\left| \int_{\overline{z_1 z_0}} f^{i-1}\, dz \right| \geq \left| f_{(Y_0,T_i)}(z_0) \int_{\overline{z_1 z_0}} dz \right| - \left| \int_{\overline{z_1 z_0}} \left( f_{(Y_0,T_i)}(z_0)-f_{(Y_0,T_i)}(z) \right)\, dz \right|
\\
- \left| \int_{\overline{z_1 z_0}} \left( f_{(Y_0,T_i)}(z)-f^{i-1}(z) \right)\, dz \right| \geq \frac{\xi}2 \left( \left| f_{(Y_0,T_i)}(z_0) \right| -\ep_1-\ep_1 \right) \geq
\end{multline*}
\[
\frac{\xi}2 ( \min_{Q_i} \{|f_{(Y_0,T_i)}|\}-2\ep_1 )=\frac{\xi}4 \min_{Q_i}\{|f_{(Y_0,T_i)}|\}\;.
\]
Then, joining this computation with \eqref{puti} and taking (f5.$i$) into account, we obtain that
\[
\normaR{(Y_n(z_0)-Y_n(z_1))_{(*,T_i)}} \geq \frac12 ( \tau_i \frac{\xi}4 \min_{Q_i}\{|f_{(Y_0,T_i)}|\}-1)-2\ep_0 > 2(\Lambda+1-\ep_0)\;.
\]
Therefore, there exists $\a\in\{0,1\}$ such that
\begin{equation}\label{BBB}
\normaR{(Y_n(z_\a)-X_n(z_0))_{(*,T_i)}} > \Lambda\;.
\end{equation}

On the other hand,
\begin{multline}\label{DDD}
|(Y_n(z_\a)-X_n(z_0))_{(3,T_i)}|\leq \normaR{Y_n(z_\a)-Y_i(z_\a)}
\\
+ |(Y_i(z_\a)-Y_{i-1}(z_\a))_{(3,T_i)}| +\normaR{Y_{i-1}(z_\a)-X_n(z_\a)} + \normaR{X_n(z_\a)-X_n(z_0)}< 3\ep_0\;,
\end{multline}
where we have used (f2.$k$), $k\neq i,$ (f1.$i$) and (e8). Hence, using \eqref{BBB} and \eqref{DDD} we conclude that we can write $Y_n(z_\a)=X_n(z_0)+u+\l v,$ where $z_0\in Q_i,$ $u\in\{x_1^2+x_2^2+x_3^2\leq 1\},$ $\l > \Lambda$ and $v\in \Pi_i$ with $\normaR{v}=1.$ Therefore, \eqref{martilla} guarantees that $Y_n(z_\a)\in \L^3\setminus \overline{B(r_2)}.$
\end{proof}

From \eqref{20}, it is clear that $Y(P)\subset B(r_2).$ Then, the existence of a polygon $Q$ satisfying items (I) and (IV) is a direct consequence of Claim \ref{Af-2}.
\\

{\bf Item (V):} Again, as a previous step, we consider the following statement. Its proof is elemental, we leave the details to the reader.

\begin{claim}\label{fuera}
Consider $z_0\in E(r_2)\setminus \overline{B(r_1)}$ and $T$ the tangent plane to $b(r_2)$ at the point $\mathcal{P}_H^{r_2}(z_0).$ Let $T_0$ be the parallel plane to $T$ passing through $z_0.$ Then,
\[
T_0\subset\L^3\setminus B (r_1-1)\;.
\]
\end{claim}

Now, we are proving that item (V) holds. Given $z\in \intc Q\setminus \intc P,$ there are five possible situations for the point $z$ (recall that $Q_i^\xi\cap D(p_j,\de)=\emptyset,$ $\forall j\notin \{i,i+1\}$).
\\

\noindent{\bf Case 1)} Assume $z\notin (\cup_{k=1}^n D(p_k,\de))\cup (\cup_{k=1}^n Q_k^\xi).$ In this case we can make use of properties (d1.$k$) and (f2.$k$), $k=1\ldots,n,$ to conclude that
\[
\normaR{Y(z)-X(z)}< 2\ep_0\;,
\]
so, if $\ep_0$ is small enough, we can finish using \eqref{6}.
\\

\noindent{\bf Case 2)} Suppose $z\in D(p_i,\de)\setminus \cup_{k=1}^n Q_k^\xi.$ In this case, we use (f2.$k$), (d1.$k$), $k=1,\ldots,n,$ \eqref{7}, (b8.$i$) and the fact that $S_i$ is an $\R^3$-orthonormal basis to obtain
\begin{multline*}
\escproR{Y_n(z)-X(p_i),e^i_2}= \escproR{Y_n(z)-Y_0(z),e^i_2}+ \escproR{X_n(z)-X_i(z),e^i_2}
\\
+ \escproR{X_i(z)-X_i(q_i),e^i_2} + \escproR{X_i(q_i)-X(q_i),e^i_2} + \escproR{X(q_i)-X(p_i),e^i_2}>
\end{multline*}
\[
\escproR{X_i(z)-X_i(q_i),e^i_2}-4\ep_0> \frac12 |f_{(X,S_i)}(p_i)|\int_{\overline{q_i z}}\frac{k_i\, dw}{w-p_i}-8\ep_0>-8\ep_0\;.
\]

In the same way, but using (d2.$i$) instead of (b8.$i$) we conclude that
\[
\escproR{Y_n(z)-X(p_i),e^i_3} > -4\ep_0\;.
\]

Again, if $\ep_0$ is sufficiently small, we can finish the proof taking into account the above inequalities and \eqref{6}.
\begin{figure}[htbp]
    \begin{center}
        \includegraphics[width=0.60\textwidth]{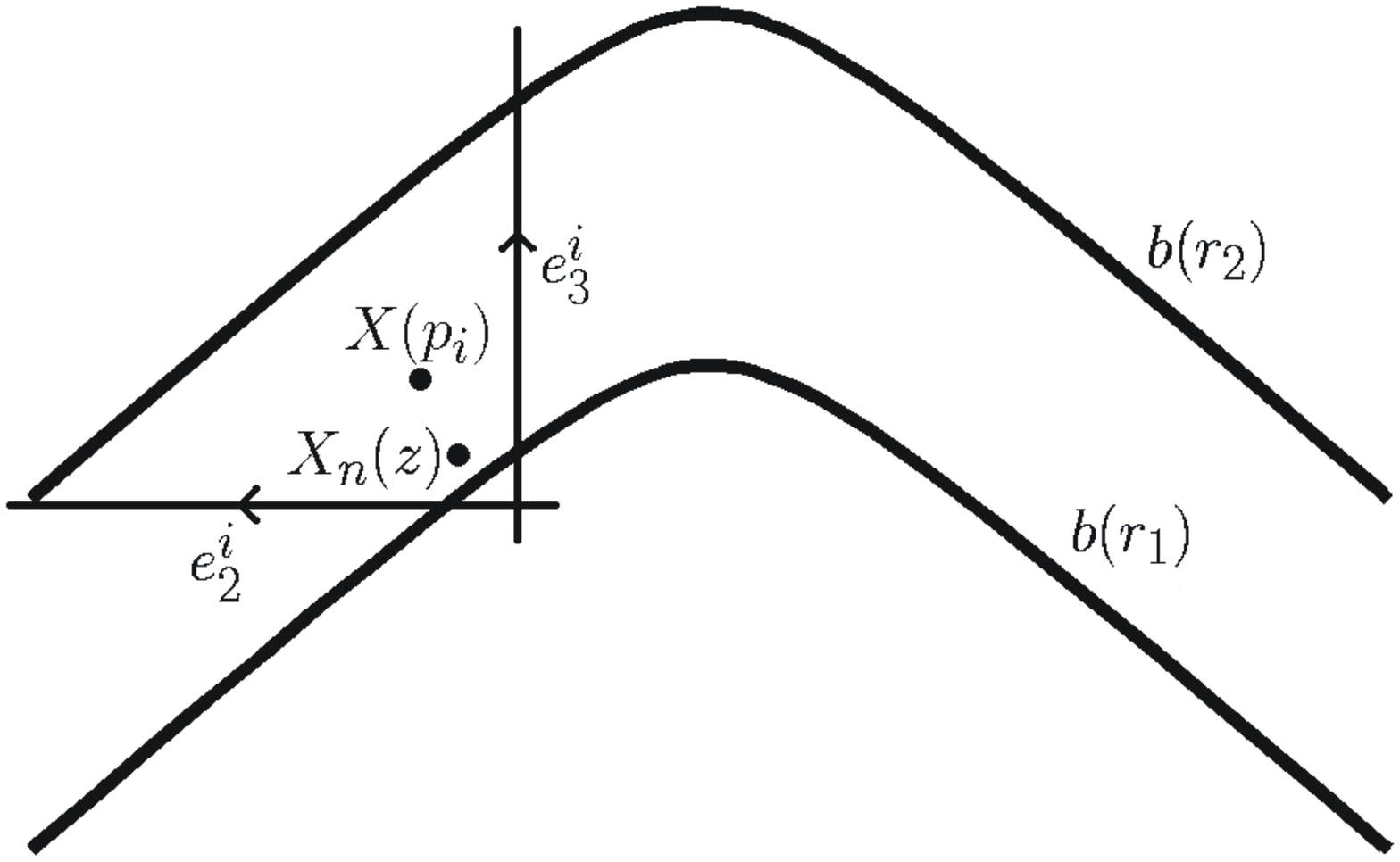}
    \end{center}
    \caption{A possible place for the point $X_n(z)$ in Case 2).} \label{fig:fuera1}
\end{figure}
\\

\noindent{\bf Case 3)} Assume $z\in D(p_i,\de)\cap Q_i^\xi.$ Following the arguments of the above case, we can obtain
\begin{equation}\label{lanj}
\escproR{X_n(z)-X(p_i),e^i_2}>-7\ep_0\;,
\end{equation}
and
\begin{equation}\label{aron}
\escproR{X_n(z)-X(p_i),e^i_3}>-3\ep_0\;.
\end{equation}

On the other hand, making use of (f2.$k$), $k=1,\ldots,n,$ and (f1.$i$) we know that
\begin{equation}\label{lapiz}
(Y_n(z)-X_n(z))_{(3,T_i)}>-\ep_0\;.
\end{equation}

Now, label $d_i:=X_n(a_i)-3\mu\, \mathcal{N}_H^{r_2}(X(p_i)).$ Then, (d5.$i$) implies
\begin{equation}\label{raton}
\normaR{X(p_i)-d_i}<14\ep_0\;.
\end{equation}
Moreover, $\mathcal{P}_H^{r_2}(d_i)=\mathcal{P}_H^{r_2}(X_n(a_i)),$ and taking into account \eqref{6} and \eqref{raton}, if $\ep_0$ is small enough, we have $d_i\in E(r_2)\setminus \overline{B(r_1)}.$ Let $\Pi$ be the tangent plane to $b(r_2)$ at the point $\mathcal{P}_H^{r_2}(X_n(a_i)).$ Given $x\in \L^3,$ denote by $\Pi_x$ the parallel plane to $\Pi$ passing through $x.$ Then, we can apply Claim \ref{fuera} to the point $d_i$ obtaining that
\[
\Pi_{d_i} \subset\L^3\setminus B(r_1-1)\;.
\]
Therefore, taking \eqref{raton} into account we conclude that
\[
\Pi_{X(p_i)} \subset\L^3\setminus B\big(r_1-1-\frac{b_2}3\big)\;,
\]
where $\ep_0$ must be chosen small enough. Hence, using \eqref{lanj} and \eqref{aron} one has (for $\ep_0$ sufficiently small)
\[
\Pi_{X_n(z)} \subset\L^3\setminus B\big(r_1-1-\frac{b_2}3-\frac{b_2}3\big)\;.
\]
\begin{figure}[htbp]
    \begin{center}
        \includegraphics[width=1\textwidth]{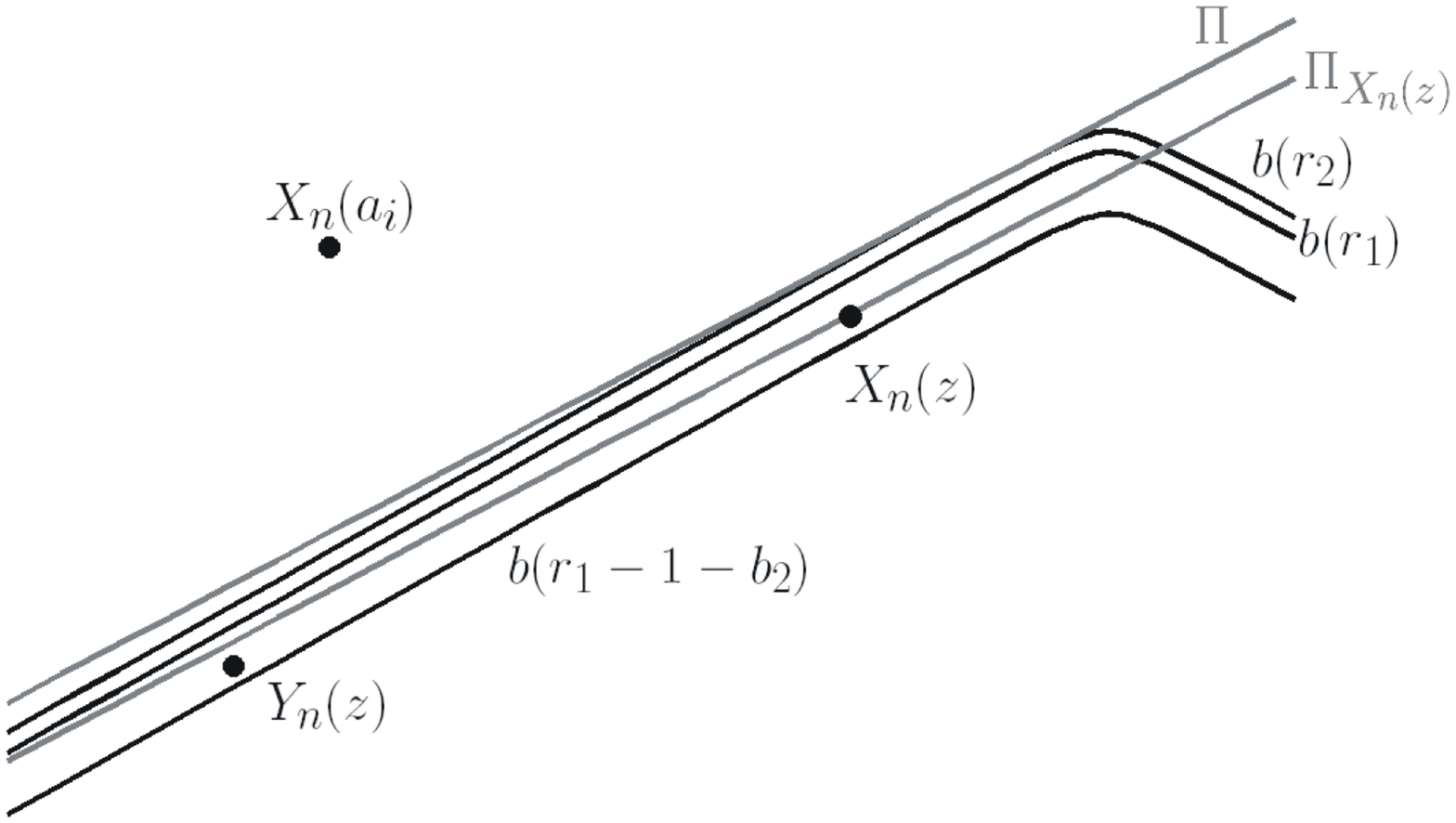}
    \end{center}
    \caption{A possible place for the point $Y_n(z)$ in Case 3).} \label{fig:fuera2}
\end{figure}
Finally, the above equation and \eqref{lapiz} guarantee that $Y_n(z)\in \L^3\setminus B(r_1-1-b_2),$ where again we have to take $\ep_0$ small enough.
\\

\noindent{\bf Case 4)} Assume $z\in D(p_{i+1},\de)\cap Q_i^\xi.$ In this case, taking also \eqref{8} into account we can obtain
\[
\escproR{X_n(z)-X(p_i),e^i_2}>-7\ep_0-\frac{\ep_0}{3\mu}\;,\quad \escproR{X_n(z)-X(p_i),e^i_3}>-3\ep_0-\frac{\ep_0}{3\mu}\;.
\]
Then, we finish reasoning as in the former case.
\\

\noindent{\bf Case 5)} Finally, assume $z\in Q_i^\xi\setminus \cup_{k=1}^n D(p_k,\de).$ Now, we can apply (d1.$k$), $k=1,\ldots,n,$ and \eqref{7} to obtain
\[
\normaR{X_n(z)-X(p_i)}<2\ep_0\;.
\]
Again, we conclude the proof reasoning as in case 3).
\\

This last case concludes the proof of item (V) and completes the proof of Lemma \ref{lema}.

\end{document}